\documentclass[11pt]{article}

\usepackage{amsmath,amssymb,latexsym}
\usepackage{mathrsfs}
\usepackage{verbatim}

\newtheorem{defn}{Definition}[section]

\newtheorem{ex}[defn]{Example}
\newtheorem{thm}[defn]{Theorem}
\newtheorem{prop}[defn]{Proposition}
\newtheorem{cor}[defn]{Corollary}

\newcommand{\h}{{H}}

\newcommand{\mn}{\mathbb N}
\newcommand{\mr}{\mathbb R}

\newcommand{\nulel}{{\bf 0}}

\newcommand{\seqgr[1]}{\{#1_i\}_{i=1}^\infty }
\newcommand{\sech[1]}{{\cap_{s\in\mn_0}{#1}_s}}

\newcommand{\snorm[1]}{{\|\hspace{-.01in}|{#1}\|\hspace{-.015in}|}}

\def\be{\begin{equation}}
\def\ee{\end{equation}}
\def\newin {\kern-0.24em\in\kern-0.19em}
\def\newsubset {\kern-0.2em\subset\kern-0.2em}
\def\vn{\vspace{.1in}\noindent}
\def\v{\vspace{.1in}}
\def\sumii{\sum_{i=1}^{\infty}}
\def\bp{\noindent{\bf Proof. \ }}
\def\ep{\noindent{$\Box$}}
\def\<{\langle}
\def\>{\rangle}

\begin{document}

\title{On sequence spaces for Fr\'echet frames}

\author{Stevan Pilipovi\'c\thanks{This research was supported by the Ministry of Science of Serbia,
Project 144016} \ and Diana T. Stoeva\thanks{This research was supported by DAAD.}}

\maketitle

\begin{abstract}
We analyze the construction of a sequence space $\widetilde{\Theta}$, resp. a sequence of sequence spaces, in order to have $\seqgr[g]$ as a $\widetilde{\Theta}$-frame or Banach frame for a Banach space $X$, resp. pre-$F$-frame or $F$-frame for a Fr\'echet space $X_F=\cap_{s\in\mn_0} X_s$, where $\{X_s\}_{s\in\mn_0}$ is a  sequence of Banach spaces.
\end{abstract}

\section{Introduction}

This paper is closely connected to \cite{PS} and we refer to \cite{PS} for the background material. In order to keep the information about the sources for our investigations, we quote the same literature as in \cite{PS}. 

Let $X$ be a Banach space (resp. $\{X_s\}_{s\in\mn_0}$ be a sequence of Banach spaces), $\Theta$ be a $BK$-space and $\seqgr[g]$ be a $\Theta$-Bessel sequence for $X$ (resp. for $X_0$). We have investigated in \cite{PS} constructions of a $CB$-space $\widetilde{\Theta}$ 
(resp. a sequence $\{\Theta_s\}_{s\in\mn_0}$ of $CB$-spaces), 
so that given $\Theta$-Bessel sequence $\seqgr[g]$ is a $\widetilde{\Theta}$-frame or Banach frame for $X$ with respect to $\widetilde{\Theta}$ (resp. pre-$F$-frame or $F$-frame for $X_F=\cap_{s\in\mn_0} X_s$ with respect to $\Theta_F=\cap_{s\in\mn_0} \Theta_s$).

In this paper we observe that one can actually start with a
sequence $\seqgr[g]$, without considering a sequence space $\Theta$ and without the assumption for $\seqgr[g]$ to be a $\Theta$-Bessel sequence.
Our motivation comes from some sequences $\seqgr[g]$ which are not Bessel sequences (and thus \cite[Theorem 4.6]{PS} do not apply to them) but they  
give rise to series expansions (see $\seqgr[g^1]$ and $\seqgr[g^2]$ in Example \ref{ex3}).

 Let $\seqgr[g]\in (X^*)^\mn$ be given and let there exist $\seqgr[f]\in X^\mn\setminus \{\nulel\}$ such that 
 the following series expansion in $X$ holds  
 \begin{equation}\label{se}
 f=\sumii g_i(f) f_i, \ f\in X.
 \end{equation}
 Validity of (\ref{se}) does not imply that $\seqgr[g]$ is a Banach frame for $X$
with respect to the given sequence space in advance, as one can see in the following examples.
\begin{ex}\label{ex3}
Let $\seqgr[e]$ be an orthonormal basis for the Hilbert space $\h$. Consider the sequences

1. $\seqgr[g^1]:=\{e_1, e_1, 2e_2, 3e_3, 4e_4, \ldots\}$,

2. $\seqgr[g^2]:=\{e_1, e_2, e_1, e_3, e_1, e_4, e_1, e_5, \ldots\}$.

\noindent Clearly, none of these sequences is a Hilbert frame for $\h$, i.e. none of these sequences is a Banach frame for $\h$ with respect to $\ell^2$.
However, series expansions in $\h$ in the form (\ref{se}) exist via the sequences

1. $\seqgr[f^1]:=\{e_1, 0, \frac{1}{2}e_2, \frac{1}{3}e_3,\frac{1}{4}e_4,\ldots\}$,

2. $\seqgr[f^2]:=\{e_1, e_2, 0, e_3, 0, e_4, 0, e_5, \ldots\}$,

\noindent respectively.
\end{ex}
Validity of (\ref{se}) implies that $\seqgr[g]$ is a Banach frame for $X$ with respect to the $CB$-space $\Theta_f$, defined by 
$$\Theta_f:=\{ \seqgr[c]\ : \ \sumii c_i f_i \ \mbox{converges in $X$}\},$$ 
$$\|\seqgr[c]\|_{\Theta_f}:=\sup_n\|\sum_{i=1}^n c_i f_i \|_X,$$ see \cite{CCS}.
Therefore, the sequences $\seqgr[g^1]$ and $\seqgr[g^2]$ from Example \ref{ex3} are Banach frames for $\h$ with respect to the corresponding sequence spaces $\Theta_{f^1}$, $\Theta_{f^2}$, respectively.
The definition of the space $\Theta_f$ involves the use of $\seqgr[f]$. 

Our aim in this paper is to find properties on $\seqgr[g]$ so that the sequence space $\widetilde{\Theta}$, which will be investigated below (see 
(\ref{thetasnorm})), is a $CB$-space and $\seqgr[g]$ is a Banach frame for $X$ with respect to $\widetilde{\Theta}$. In assertions (\ref{constructone2})-(\ref{hilbcornew}), we determine such properties on $\seqgr[g]$. 
It is not completely answered yet whether in general the expansion property (\ref{se}) with $\seqgr[g]$ and corresponding $\seqgr[f]$ implies that 
$\widetilde{\Theta}$ is $CB$ and that $\seqgr[g]$ is a Banach frame for $X$ with respect to $\widetilde{\Theta}$.

\v
Now we will describe the results of the present paper.
In Theorem \ref{constructone2} we construct $\widetilde{\Theta}$ such that $\seqgr[g]$ is a $\widetilde{\Theta}$-Bessel sequence for $X$. Actually, in Theorem \ref{constructone2} we consider several conditions ($\mathcal{A}_1$)-($\mathcal{A}_3$) in order to 
characterize  
$\widetilde{\Theta}$-Bessel sequence for $X$ or Banach frame for $X$ with respect to $\widetilde{\Theta}$. 
Further on, we construct a class of Banach frames and pre-$F$-frames started with a Hilbert space and its orthonormal basis. Proposition \ref{primer1new} deals with possible designs of a Banach frame which shows that conditions ($\mathcal{A}_1$)-($\mathcal{A}_3$) are intrinsically related to our construction of $\widetilde{\Theta}$. 
 Theorem \ref{constructthetas2} and Corollary \ref{c1} deal with the construction of sequence spaces $\Theta_s$ (as $\widetilde{\Theta}$, but with $X_s$, $s\in\mn_0$) so that one obtains a pre-$F$-frame or $F$-frame for $X_F=\cap_{s\in\mn_0} X_s$. 
Proposition \ref{hilbcornew} deals with a sequence of Hilbert spaces and a sequence $\seqgr[g]$ designed similarly as in Proposition \ref{primer1new}, in order to analyze conditions ($\mathcal{A}_1^s$)-($\mathcal{A}_3^s$) and the sequence spaces $\Theta_s$, $s\in\mn_0$, so that $\seqgr[g]$ is a pre-$F$-frame or $F$-frame for $X_F$. 

\section{Pre-$F$- and $F$-frames}

We will use notation and notions as in \cite{PS}.
$(X, \|\cdot \|)$ is a Banach space and
$(X^*, \|\cdot \|^*)$ is its dual, $( \Theta, \snorm[\cdot]) $ is
a Banach sequence space and  $( \Theta^*, \snorm[\cdot]^*) $ is
the dual of $\Theta$. Recall that $\Theta$ is called {\it solid}
if the conditions $\{c_i\}_{i=1}^\infty\in \Theta$ and $|d_i| \leq
|c_i|,$ $i\newin\mn$, imply that $\{d_i\}_{i=1}^\infty \in \Theta$
and $\snorm[\{d_i\}_{i=1}^\infty] \leq
\snorm[\{c_i\}_{i=1}^\infty].$ If the coordinate functionals on
$\Theta$ are continuous, then $\Theta$ is called a {\it
$BK$-space}. A $BK$-space, for which the canonical vectors form a
Schauder basis, is called a {\it $CB$-space}. 

Let $\{Y_s, | \cdot |_s\}_{s\in\mn_0}$ be a sequence of separable
Banach spaces such that \be \label{fx1} \{\nulel\} \neq
\sech[Y]\subseteq \ldots \subseteq Y_2 \subseteq Y_1 \subseteq Y_0
\ee \be  \label{fx2} |\cdot|_0\leq | \cdot |_1\leq | \cdot |_2\leq
\ldots \ee \be \label{fx3} Y_F :=\sech[Y] \;\; \mbox{is dense in}
\;\;\; Y_s, \;\;\; \mbox{for every} \;\;\; s\in\mn_0. \ee

Then $Y_F$ is a Fr\'echet space
with the sequence of norms $ | \cdot |_s, $ $ s\in\mn_0.$
 We use the above sequences in two cases:

1. $Y_s=X_s$ with norm $\|\cdot\|_s, s\in\mn_0;$

2. $Y_s=\Theta_s$ with norm $\snorm[\cdot]_s, s\in\mn_0$.

Let $\{X_s, \|\cdot\|_s\}_{s\in\mn_0}$ and $\{\Theta_s,
\snorm[\cdot]_s\}_{s\in\mn_0}$ be sequences of Banach spaces,
which satisfy (\ref{fx1})-(\ref{fx3}). For fixed $s\newin\mn_0$,
an operator $V:\Theta_F\to X_F$ is called $s$-bounded, if
there exists a constant $K_s>0$ such that $\|V\seqgr[c]\|_s\leq
K_s \snorm[\{c_i\}_{i=1}^\infty]_s$ for all $\seqgr[c]\newin
\Theta_F$. If $V$ is $s$-bounded for every $s\newin\mn_0$, then
$V$ is called $F$-bounded.

\begin{defn}\label{fframe}
Let $\{X_s, \|\cdot\|_s\}_{s\in\mn_0}$ be a sequence of Banach spaces
which satisfies (\ref{fx1})-(\ref{fx3}), and let $\{\Theta_s,
\snorm[\cdot]_s\}_{s\in\mn_0}$ be a sequence of $BK$-spaces
which satisfies (\ref{fx1})-(\ref{fx3}). A sequence $\seqgr[g]\newin
({X_F^*})^\mn$ is called a pre-$F$-frame for $X_F$ with respect to
$\Theta_F$ if for every $s\in\mn_0$ there exist constants $0<A_s\leq
B_s<\infty$ such that
\begin{equation}\label{fframestar}
\{g_i(f)\}_{i=1}^\infty\in\Theta_F, \ f\newin X_F,
\end{equation}
\begin{equation}\label{fframetwostar}
A_s \|f\|_s\leq \snorm[\{g_i(f)\}_{i=1}^\infty]_{s}\leq B_s\|f\|_s, \ f\newin X_F.
\end{equation}
The constants $B_s$ (resp. $A_s$),
$s\in\mn_0$, are called upper (resp. lower) bounds for
$\seqgr[g]$. The pre-$F$-frame is called {\it tight}, if $A_s=B_s,
s\newin\mn_0$.

Moreover, if there exists an
$F$-bounded operator
$V:\Theta_F\rightarrow X_F$ so that $V(\{g_i(f)\}_{i=1}^\infty)=f$
for all $f\newin X_F,$ then a pre-$F$-frame $\seqgr[g]$ is called
an $F$-frame (Fr\'echet frame) for $X_F$ with respect to $\Theta_F$ and $V$ is
called an $F$-frame operator for $\seqgr[g]$.

When (\ref{fframestar}) and at least the upper inequality in
(\ref{fframetwostar}) hold,
then $\seqgr[g]$ is
called an $F$-Bessel sequence for $X_F$ with respect to $\Theta_F$
with bounds $B_s$, $s\newin\mn_0$.
\end{defn}

\begin{thm} \label{constructone} {\rm \cite{PS}}
Let $\Theta\neq\{\nulel\}$ be a solid $BK$-space,
$X\neq\{\nulel\}$ be a reflexive Banach space and $\seqgr[g]\newin
(X^*)^\mn$ be a $\Theta$-Bessel sequence for $X$ with  bound
$B\leq 1$ such that $0<\|g_i\|\leq 1$, $i\newin\mn$. For every
$c=\seqgr[c]\newin \Theta$, denote
\begin{equation} \label{setm}
M^c:=\{f\in X \ : \ |c_i|\leq |g_i(f)|, \ i\in \mn \}
\end{equation} and define
\begin{equation} \label{setm2}
 \widetilde{\Theta} := \left\{ c\in\Theta \ : \ M^{c}\neq
\varnothing\right\}, \ \, \snorm[c]_{\widetilde{\Theta}}
:=\inf\left\{ \|f\| : f\in M^{c}\right\}.
\end{equation}

Consider the conditions:

\vspace{.1in} \noindent $(\mathcal{A}_1): \ (\forall
\,c\in\widetilde{\Theta})\ (\forall \, d\in \widetilde{\Theta}) \
(\forall\, f\in M^c) \ (\forall\, h\in M^d) \Rightarrow $

\vspace{.1in} \hspace{1.7in} $ (\exists \ r\in M^{c+d})\ \
(\|r\|\leq \|f\|+\|h\|).$

\vn $(\mathcal{A}_2): \ (\forall c\in \widetilde{\Theta}) \
(\forall\, \varepsilon>0) \ (\exists \, k\in\mn) \  (\exists
\,f\in M^{c^{(k)}})  \ ( \|f\|<\varepsilon).$

Assume that $(\mathcal{A}_1)$ is valid. Then the following holds.
 \begin{itemize}
\item[{\rm(a)}]
$\widetilde{\Theta}$ is a solid $BK$-space
with
$\snorm[\cdot]_{\Theta}\leq\snorm[\cdot]_{\widetilde{\Theta}}$ and
$\{g_i\}_{i=1}^\infty$ is a $\widetilde{\Theta}$-Bessel sequence
for $X$ with bound $\widetilde{B}=1$.
\item[{\rm(b)}] If $\seqgr[g]$ is a $\Theta$-frame for
$X$, then $\{g_i\}_{i=1}^\infty$ is a $\widetilde{\Theta}$-frame
for $X$.
\item[{\rm(c)}] If $\seqgr[g]$ is a Banach frame for
$X$ with respect to $\Theta$, then $\{g_i\}_{i=1}^\infty$ is a
Banach frame for $X$ with respect to $\widetilde{\Theta}$.
\item[{\rm(d)}] $\widetilde{\Theta}$ is a $CB$-space if and only if
$(\mathcal{A}_2)$ holds.
\end{itemize}
\end{thm}

\begin{thm}\label{constructthetas} {\rm \cite{PS}}
  Let $\Theta\neq\{\nulel\}$ be a solid
$BK$-space and $\{X_s\}_{s\in\mn_0}$ be a sequence of reflexive
Banach spaces which satisfies (\ref{fx1})-(\ref{fx3}). Let
$\seqgr[g]\newin (X_0^*)^\mn$ be a $\Theta$-Bessel sequence for
$X_0$ with bound $B\leq 1$ such that $0<\|g_i\|\leq 1$,
$i\newin\mn$. For every $s\newin\mn_0$ and every
$c=\seqgr[c]\newin \Theta$, denote
\begin{equation}\label{mx2}
M^{c}_s:=\{f\in X_s \ : \ |c_i|\leq |g_i(f)|, \ i\in \mn \}
\end{equation}
and define
 \begin{equation}\label{ts2}
 {\Theta}_s := \left\{ c\in\Theta \ : \ M^c_s\neq
\varnothing\right\}, \ \, \snorm[c]_{s} :=\inf\left\{ \|f\|_s :
f\in M^c_s\right\}.
 \end{equation}
Consider the following conditions:

\vspace{.1in}\noindent $(\mathcal{A}_1^s):$ \ \ \ \ \
 $(\forall c\in\Theta_s)\ \ (\forall d\in \Theta_s)\ \ (\forall f\in
M^c_s) \ \ (\forall h\in M^d_s) \ \Rightarrow $

\vspace{.1in} \centerline{$(\exists \,r\in M^{c+d}_s)\ \
(\|r\|_s\leq \|f\|_s+\|h\|_s).$}

\vspace{.1in}\noindent $(\mathcal{A}_2^s): \
 \ (\forall c\newin \Theta_s) \  (\forall\,
\varepsilon>0)\ \ (\exists \, k\newin\mn) \  \, (\exists \,f\newin
M^{c^{(k)}}_s) \, \ (\|f\|_s<\varepsilon).$

\vspace{.1in}\noindent $(\mathcal{A}_3^s):$ \ There exists
$A_s\in(0,1]$ such that for every $f\in X_s$ one has
\begin{equation}\label{fj1}
\widetilde{f}\in M_s^{\{g_i(f)\}_{i=1}^\infty} \Rightarrow
A_s\|f\|_s\leq \|\widetilde{f}\|_s.
\end{equation}

Assume that $(\mathcal{A}_1^s)$ holds for every $s\in\mn_0$.
Then the following holds.

\begin{itemize}
\item[$(\mathcal{P}_1)$] $\{\Theta_s\}_{s\in\mn_0}$ is a sequence of solid
$BK$-spaces with the properties (\ref{fx1})-(\ref{fx2}) such that
$\{g_i|_{X_s}\}_{i=1}^\infty$ is a $\Theta_s$-Bessel sequence for
$X_s$ with bound $B_s=1$, $s\newin\mn_0$.
\item[$(\mathcal{P}_2)$] For any $s\newin\mn$,
 $\{g_i|_{X_s}\}_{i=1}^\infty$ is a
$\Theta_s$-frame for $X_s$ if and only if $(\mathcal{A}_3^s)$
holds.
If $(\mathcal{A}_3^s)$ holds with $A_s=1$, then
$\{g_i\vert_{X_s}\}$ is a tight $\Theta_s$-frame for $X_s$.
\item[$(\mathcal{P}_3)$] For any $s\newin\mn$,
$\Theta_s$ is a $CB$-space if and only if $(\mathcal{A}_2^s)$
holds.
\end{itemize}

\end{thm}

\section{Construction of a sequence space}

We construct a sequence space $\widetilde{\Theta}$ and through the properties $(\mathcal{A}_1)$-$(\mathcal{A}_3)$ we analyze frame-properties of a given sequence $\seqgr[g]$ and given $X$.

\begin{thm} \label{constructone2}
Let $(X, \|\cdot\|)\neq\{\nulel\}$ be a reflexive Banach space
and $\seqgr[g]\newin (X^*)^\mn\setminus \{\nulel\}$.
For every scalar sequence $c=\seqgr[c]$, let $M^c$ be given by (\ref{setm}) and
define
 \begin{equation}\label{thetasnorm}
 \widetilde{\Theta} := \left\{ c=\seqgr[c] \ : \ M^{c}\neq
\varnothing\right\}, \ \, \snorm[c]_{\widetilde{\Theta}}
:=\inf\left\{ \|f\| : f\in M^{c}\right\}.
 \end{equation}
Consider conditions $(\mathcal{A}_1)$, $(\mathcal{A}_2)$ and  

{\rm ($\mathcal{A}_3$)} : there exists
$A\in(0,1]$ such that for every $f\in X$ one has
\begin{equation*}\label{fj1new} 
\widetilde{f}\in M^{\{g_i(f)\}_{i=1}^\infty} \Rightarrow
A\|f\|\leq \|\widetilde{f}\|.
\end{equation*}
Assume that $(\mathcal{A}_1)$ is fulfilled.
 Then the following holds.

\vspace{.05in} {\rm(a)} $\widetilde{\Theta}$ is a solid $BK$-space
 and $\{g_i\}_{i=1}^\infty$ is a $\widetilde{\Theta}$-Bessel sequence for $X$
with bound $\widetilde{B}=1$.

\vspace{.05in} {\rm(b)} $\{g_i\}_{i=1}^\infty$ is a
$\widetilde{\Theta}$-frame for $X$ if and only if ($\mathcal{A}_3$) holds
if and only if there exists a solid $BK$-space $\Theta\supseteq \widetilde{\Theta}$
such that $\seqgr[g]$ is a $\Theta$-frame for $X$.

{\rm(c)} $\seqgr[g]$ is a Banach frame for $X$ with respect to $\widetilde{\Theta}$ if and only if there exists a solid $BK$-space
$\Theta\supseteq \widetilde{\Theta}$ such that $\seqgr[g]$ is a Banach frame for $X$ with respect to $\Theta$.

{\rm(d)} $\widetilde{\Theta}$ is a $CB$-space if and only if
$(\mathcal{A}_2)$ holds.
\end{thm}

\bp (a) In the same way as in \cite[Theorem 4.6]{PS}, using ($\mathcal{A}_1$) it follows that $\widetilde{\Theta}$ is a linear space and
$\snorm[\cdot]_{\widetilde{\Theta}}$
 is a norm, only the proof that $\snorm[c]_{\widetilde{\Theta}}=0$ implies $c_i=0$, $i\newin\mn$, is different. Let $\snorm[c]_{\widetilde{\Theta}}=0$. Fix $i\in\mn$. For every
$\varepsilon>0$, there exists $f_\varepsilon\in M^c$ such that $\|f_\varepsilon\|<\varepsilon/{\|g_i\|}$ and hence $|c_i|\leq
|g_i(f_\varepsilon)|<\varepsilon$ which implies that $c_i=0$.

For solidity-property, let $c=\seqgr[c]\in
\widetilde{\Theta}$ and $d=\seqgr[d]$ be such that $|d_i|\leq
|c_i|$, $i\newin\mn$. Since $M^{d}\supseteq M^{c}\neq \varnothing$, it follows that
$\seqgr[d]\in\widetilde{\Theta}$ and
$\snorm[d]_{\widetilde{\Theta}}\leq\snorm[c]_{\widetilde{\Theta}}$.

For the completeness of $\widetilde{\Theta}$, first note that the $i$-th coordinate functional on $\widetilde{\Theta}$ is continuous. Indeed, fix $i\in\mn$. If
$c=\{c_j\}_{j=1}^\infty\in \widetilde{\Theta}$, then for every $f\in M^c$ we have $|c_i|\leq |g_i(f)|\leq \|g_i\| \|f\|$ and hence
\begin{equation}\label{coordfunct} |c_i|\leq \|g_i\| \inf\left\{ \|f\| : f\in M^{c}\right\}=\|g_i\| \snorm[c]_{\widetilde{\Theta}}. \end{equation}
 Let
now $c^\nu=\{c_i^\nu\}_{i=1}^\infty, \nu\in\mn,$ be a Cauchy
sequence in $\widetilde{\Theta}$. Fix arbitrary $\varepsilon
>0$. There exists $\nu_0(\varepsilon)$ such that for
every $\mu,\nu\in\mn$, $\mu\geq\nu_0$, $\nu\geq \nu_0$, there
exists $f^{\mu,\nu}\in X$, such that
\begin{equation*}\label{compl1anew}
 \|f^{\mu,\nu}\|<\varepsilon \ \mbox{and} \ |c_i^\mu -
c_i^\nu|\leq  |g_i(f^{\mu,\nu})|, \
 i\in\mn.
\end{equation*}
By (\ref{coordfunct}), for every $i\in\mn$ the sequence $c_i^\nu, \nu\newin\mn$, is a Cauchy sequence in $\mr$ and hence it converges to some
number $c_i$ when $\nu\to\infty$. Denote $c:=\seqgr[c]$. Fix $\nu\geq \nu_0$. Now, in the same way as in \cite[Theorem 4.6]{PS}, there
exists $F^\nu\in M^{c-c^\nu}$ with $\|F^\nu\|\leq\varepsilon$, which implies that $c\in\widetilde{\Theta}$ and
$$\snorm[\{c_i-c_i^\nu\}_{i=1}^\infty]_{\widetilde{\Theta}}\leq  \|F^\nu\|\leq \varepsilon, \ \nu\geq \nu_0.$$
For the sake of completeness we add brief sketch of the arguments.
Since $\|f^{\mu,\nu}\|<\varepsilon $ for every
$\mu\geq\nu_0$, there exists a convergent subsequence
$\{\|f^{\mu_k,\nu}\|\}_{k=1}^\infty$; denote its limit by $a^\nu$.
Since $X$ is reflexive and the sequence
$\{f^{\mu_k,\nu}\}_{k=1}^\infty$ is norm-bounded, by
\cite[Corollary 1.6.4]{AK} there exists a subsequence
$\{f^{\mu_{k_n},\nu}\}_{n=1}^\infty$ which converges weakly to
some element $F^\nu\newin X$. Therefore,
$$\|F^\nu\|\leq
\lim \inf
\|f^{\mu_{k_n},\nu}\|=\lim_{n\to\infty}\|f^{\mu_{k_n},\nu}\|=a^\nu\leq 
\varepsilon.$$
Now, for every $i\newin\mn$, $ |c_i - c_i^\nu|\leq
|g_i(F^\nu)|$. Thus, $F^\nu$ belongs to $M^{c-c^\nu}$. 
This concludes the proof that $\widetilde{\Theta}$ is complete.

The $\widetilde{\Theta}$-Bessel property of $\seqgr[g]$ is clear by (\ref{thetasnorm}).

(b)-(c) For the first equivalence in (b), assume that ($\mathcal{A}_3$) holds. In this case we have
$$A\|f\|\leq \inf\{\|\widetilde{f}\| \ : \ \widetilde{f}\in
M^{\{g_i(f)\}_{i=1}^\infty}\}=\snorm[\{g_i(f)\}]_{\widetilde{\Theta}}, \ \forall f\in X,$$ and hence, by (a), $\seqgr[g]$ is a
$\widetilde{\Theta}$-frame for $X$.

Conversely, assume that $\{g_i\}_{i=1}^\infty$ is a
$\widetilde{\Theta}$-frame for $X$ with lower bound $A\in (0,1]$.
In this case for every $f\in X$ we have
$$A\|f\|
\leq \snorm[\{g_i(f)\}]_{\widetilde{\Theta}} \leq \|\widetilde{f}\|, \ \forall \widetilde{f}\in M^{\{g_i(f)\}_{i=1}^\infty} .$$

For the second equivalence in (b) and for (c), one of the directions is obvious, take $\Theta\equiv\widetilde{\Theta}$. For the other direction,
assume that $\Theta$ is a solid $BK$-space such that $\Theta\supseteq \widetilde{\Theta}$ and $\seqgr[g]$ is a $\Theta$-frame for $\Theta$ with upper bound $B$. Let $c\in\Theta$. For every $f\in M^c$, the solidity of $\Theta$ implies that
$\snorm[\{c_i\}_{i=1}^\infty]_\Theta \leq \snorm[\{g_i(f)\}_{i=1}^\infty]_\Theta \leq B\|f\|$. Therefore $\snorm[\{c_i\}_{i=1}^\infty]_\Theta \leq B
\snorm[\{c_i\}_{i=1}^\infty]_{\widetilde{\Theta}}$. 
The rest is similar to the proof of \cite[Theorem 4.6 (b)(c)]{PS}, but for the sake of completeness we add a proof here.
Let $A$ denote a lower $\Theta$-frame bound for the $\Theta$-frame $\seqgr[g]$. Then clearly, 
$A\|f\|
\leq B \snorm[\{g_i(f)\}_{i=1}^\infty]_{\widetilde{\Theta}},$ $f\in X$. This implies that $\seqgr[g]$ satisfies the lower
$\widetilde{\Theta}$-frame inequality.

Assume now that $\seqgr[g]$ is a Banach frame for $X$ and let $V:\Theta\to X$ denote a bounded operator such that
$V(\{g_i(f)\}_{i=1}^\infty)=f$, $f\newin X$. For every $c\in {\widetilde{\Theta}}$,
$\|V c\|\leq \|V\|\,\snorm[c]_\Theta \leq B \,\|V\|\,
\snorm[c]_{\widetilde{\Theta}},$
which implies that $V\vert_{\widetilde{\Theta}}$ is bounded on
$\widetilde{\Theta}$. 
This completes the proof.

(d) We use the same arguments as in \cite[Theorem 4.6]{PS}, but for the sake of completeness we add sketch of the proof.
First note that for any $i\newin \mn$, $g_i$ is not the null functional and hence
 the $i$-th canonical vector belongs to $\widetilde{\Theta}$.
Let $(\mathcal{A}_2)$ hold. Fix $c\in\widetilde{\Theta}$ and $\varepsilon>0$. There exists $k\in\mn$ such that
$\snorm[c^{(k)}]_{\widetilde{\Theta}}<\varepsilon$. Therefore
$\snorm[c^{(n)}]_{\widetilde{\Theta}}\leq \snorm[c^{(k)}]_{\widetilde{\Theta}}<\varepsilon$ for every
$n\geq k$, which implies that $\sum_{i=1}^n c_i e_i\to c$ in
$\widetilde{\Theta}$ as $n\to\infty.$
The converse is trivial.
\ep

 \section{Class of $\widetilde{\Theta}$-frames}
 
In this section we will consider a class of $\widetilde{\Theta}$-frames.
Actually, we generalize \cite[Proposition 4.7]{PS}, in which the sequence $\seqgr[g]$ was given by $g_1=e_1$, $g_i=e_{i-1}$, $i\in\mn, i>1$,
where $\seqgr[e]$ denoted an orthonormal basis. Now we consider the case, when every $e_i$, $i\in\mn$, can be repeted $k_i$-times.

We will use the following notation related to a sequence $c=\seqgr[c]$:
$$c^{(n)}:=\{\underbrace{0,\ldots,0}_n,c_{n+1},c_{n+2},c_{n+3},\ldots\}=c-\sum_{i=1}^n c_i e_i, \ n\in\mn,$$
 \hspace{.51in} $c^{(n)}_i$:= the $i$-th coordinate of $c^{(n)}$,
$i\in\mn$. 

\begin{prop}\label{primer1new}
Let $(X, \<\cdot,\cdot\>)$ be a Hilbert space, $\seqgr[e]$ be an
orthonormal basis for $X$ and $\Theta=\ell^2$. Let
$\seqgr[g]\newin (X^*)^\mn$ be defined by
\begin{eqnarray*}
g_i(f):=t_i\<f, e_1\>, &\ & i=1,2,\ldots, k_1, \\
g_i(f):=t_i\<f, e_j\>, &\ & i=k_{j-1}+1, k_{j-1}+2, \ldots, k_j, \ j\in\mn, j>1,
\end{eqnarray*}
where $k_j\in\mn$, $t_j\in\mr$, $t_j\neq 0$, $j\in\mn$, and let $\widetilde{\Theta}$ be defined by (\ref{setm2}).
  Then $\widetilde{\Theta}$ is a $CB$-space and $\seqgr[g]$ is a Banach frame for $X$ with respect to $\widetilde{\Theta}$.
 \end{prop}
 \bp      
Let $c=\seqgr[c]\in\widetilde{\Theta}, d=\seqgr[d]\in\widetilde{\Theta}$. Denote 
\begin{equation}\label{c1d1}
\widetilde{c}_1:=\max\left(\frac{|c_1|}{|t_1|}, \frac{|c_2|}{|t_2|}, \ldots, \frac{|c_{k_1}|}{|t_{k_1}|}\right), \  \widetilde{d}_1:=\max\left(\frac{|d_1|}{|t_1|},
\frac{|d_2|}{|t_2|}, \ldots, \frac{|d_{k_1}|}{|t_{k_1}|}\right), 
\end{equation}
\begin{equation}\label{cj}
\widetilde{c}_j:=\max\left(\frac{|c_{k_{j-1}+1}|}{|t_{k_{j-1}+1}|}, \frac{|c_{k_{j-1}+2}|}{|t_{k_{j-1}+2}|}, \ldots, \frac{|c_{k_j}|}{|t_{k_j}}|\right), \ j\in\mn, j>1,
\end{equation}
\begin{equation}\label{dj}
\widetilde{d}_j:=\max\left(\frac{|d_{k_{j-1}+1}|}{|t_{k_{j-1}+1}|},
\frac{|d_{k_{j-1}+2}|}{|t_{k_{j-1}+2}|}, \ldots, \frac{|d_{k_j}|}{|t_{k_j}|}\right), \ j\in\mn, j>1.
\end{equation}
 Fix $f\in M^c$ and $h\in M^d$. By (\ref{setm}),
 $|c_i|\leq |g_i(f)|$ and $|d_i|\leq |g_i(h)|$, $i\in \mn.$
 Therefore, 
\begin{equation} \label{nova}
\widetilde{c}_i\leq |\<f, e_i\>| \ \mbox{and} \  \widetilde{d}_i\leq |\<h, e_i\>|, \ i\in \mn.
\end{equation}
 Let us find $r^{f,h}\in M^{c+d}$ such that $\|r^{f,h}\|\leq \|f\|+\|h\|$.
Consider
$$r^{f,h}:= m_1 e_1 + m_2 e_2 + m_3 e_3 + m_4 e_4
+\ldots,$$ where
\begin{eqnarray}
m_1\!\!\!\!&=&\!\!\!\!\max \left(\frac{|c_1+d_1|}{|t_1|}, \frac{|c_2+d_2|}{|t_2|},\ldots,\frac{|c_{k_1}+d_{k_1}|}{|t_{k_1}|}\right), \label{m1}\\
m_j\!\!\!\!&=&\!\!\!\!\max \left(\frac{|c_{k_{j-1}+1}+d_{k_{j-1}+1}|}{|t_{k_{j-1}+1}|},
\frac{|c_{k_{j-1}+2}+d_{k_{j-1}+2}|}{|t_{k_{j-1}+2}|},\ldots,\frac{|c_{k_j}+d_{k_j}|}{|t_{k_j}|}\right), \label{mj}
\end{eqnarray}
for $j\in\mn, j>1$.
By (\ref{nova}), $\seqgr[\widetilde{c}]\in\ell^2$ and $\seqgr[\widetilde{d}]\in\ell^2$. Since $m_i\leq \widetilde{c}_i + \widetilde{d}_i$, $i\in\mn$, it follows that $\seqgr[m]\in\ell^2$
and hence, $r^{f,h}\in X$.
 It is clear that \begin{eqnarray*}
|c_i+d_i|&\leq& t_i m_1=t_i \<r^{f,h},e_1\> =g_i(r^{f,h}), \ \ i=1,2,\ldots, k_1,\\
|c_i+d_i|&\leq& t_i m_j=t_i \<r^{f,h},e_j\> =g_i(r^{f,h}), \ \ i=k_{j-1}+1, k_{j-1}+2,\ldots, k_j,
\end{eqnarray*}
for $j\in\mn$, $j>1$. Therefore, $r^{f,h}\newin M^{c+d}$. 
Using the facts that $m_j\leq \widetilde{c}_j + \widetilde{d}_j$, 
$\ell^2$ is solid and (\ref{nova}) holds, we obtain
that
\begin{eqnarray*}
\|r^{f,h}\| &= & \|\{m_1, m_2, m_3, \ldots \} \|_{\ell^2}\\
 &\leq & \|\{\widetilde{c}_1+\widetilde{d}_1, \widetilde{c}_2+\widetilde{d}_2, \widetilde{c}_3+\widetilde{d}_3, \ldots \} \|_{\ell^2}\\
& \leq & \|\{\widetilde{c}_1, \widetilde{c}_2,\widetilde{c}_3, \ldots \}\|_{\ell^2} +\|\{\widetilde{d}_1,\widetilde{d}_2,\widetilde{d}_3 , \ldots
\}
\|_{\ell^2}\\
& \leq & \|\{|\<f,e_i\>|\}_{i=1}^\infty\|_{\ell^2}
+\|\{|\<h,e_i\>| \}_{i=1}^\infty \|_{\ell^2}= \|f\|+\|h\|.
\end{eqnarray*}
Therefore, $(\mathcal{A}_1)$ is fulfilled.

Let us now prove that $(\mathcal{A}_3)$ is fulfilled. Consider $f\newin X$ and take arbitrary $\widetilde{f}\in M^{\{g_i(f)\}_{i=1}^\infty}$, i.e.
$|g_i(f)|\leq |g_i(\widetilde{f})|, i\in\mn$, and hence $|\<f,e_i\>|\leq |\< \widetilde{f},e_i\>|$, $i\in\mn$. Then
$$
\|f\|^2=\sumii  \,|\<f, e_i\>|^2 \leq \sumii  \,|\<\widetilde{f}, e_i\>|^2 =  \|\widetilde{f}\|_s^2.
$$
Therefore $(\mathcal{A}_3)$ holds with $A_s=1$.

Consider now $c=\seqgr[c]\in\widetilde{\Theta}$. Then $\seqgr[\widetilde{c}]\in\ell^2$. Fix $\varepsilon>0$ and find $p\newin\mn$, $p>1,$ such
that $\sum_{i=p}^\infty \widetilde{c}_i^{\,2}<\varepsilon$. Define
$$\seqgr[b]:=\{\underbrace{0, \ldots, 0}_{p-1}, \widetilde{c}_p, \widetilde{c}_{p+1},\widetilde{c}_{p+2},\ldots \}.$$
Since $\seqgr[b]\in\ell^2$, there exists $h\newin X$ such that $b_i=\<h,e_i\>$, $i\in\mn$, and $\|h\|^2=\sum_{i=1}^\infty |b_i|^2 =
\sum_{i=p}^\infty \widetilde{c}_i^2< \varepsilon.$ Take $k:=k_{p-1}$ and recall,
$$c^{(k)}=\{\underbrace{0,\ldots,0}_{k_{p-1}},c_{k_{p-1}+1},c_{k_{p-1}+2},\ldots, c_{k_p},c_{k_p+1}, \ldots\}.$$
Thus we have
\begin{eqnarray*}
|c^{(k)}_i|\!\!\!&=&\!\!\!0\leq |g_{i}(h)|, \ i=1,2,\ldots,k_{p-1}, \\
 |c^{(k)}_i|\!\!\!&\leq&\!\!\!\widetilde{c}_{p+n}=b_{p+n}=\<h,e_{p+n}\>= g_{i}(h),
 \ i=k_{p+n-1}+1, \ldots, k_{p+n}, \ n\in\mn_0,
\end{eqnarray*}
which implies that $h\in M^{c^{(k)}}$. Since $\|h\|^2<\varepsilon $, it follows that $(\mathcal{A}_2)$ is fulfilled.
Now Theorem \ref{constructone2} implies that  $\widetilde{\Theta}$ is a $CB$-space and $\seqgr[g]$ is a $\widetilde{\Theta}$-frame for $X$.

It remains to prove the existence of a Banach-frame operator. 
Let us 
denote the canonical basis for $\widetilde{\Theta}$ by $z_i$, $i\in\mn$.
Consider 
\begin{equation}\label{seqf}
\seqgr[f]:=\{\frac{1}{t_1}e_1,\underbrace{0,\ldots,0}_{k_1-1}\,, \frac{1}{t_{k_1+1}}e_2,\underbrace{0,\ldots,0}_{k_2-1}\,, \frac{1}{t_{k_2+1}}e_3,\underbrace{0,\ldots,0}_{k_3-1}\,,\ldots.\}\end{equation}
Define $V$ on $\seqgr[z]$ by $Vz_i:=f_i$, $i\in\mn$. Our aim is to prove that $V$ is bounded on $\seqgr[z]$ and then to consider the extension of $V$ on $\widetilde{\Theta}$  by linearity and continuity.
Let us first prove that $\snorm[z_i]_{\widetilde{\Theta}}=1/\|g_i\|^*$, $i\in\mn$. Define 
\begin{eqnarray}
h_i(f):=\frac{1}{t_i} e_1, &\ & i=1,2,\ldots, k_1, \label{h1}\\
h_i(f):=\frac{1}{t_i} e_j, &\ & i=k_{j-1}+1, k_{j-1}+2, \ldots, k_j, \ j\in\mn, j>1.\label{hj}
\end{eqnarray}
  Fix $i\in\mn$ and note that $h\in M^{z_i}$ if and only if $1\leq |g_i(h)|$. Thus, the elements $h\in M^{z_i}$ should satisfy the inequality $\|h\|\geq \frac{1}{\|g_i\|^*}$. 
Since $h_i\in M^{z_i}$ and $\|h_i\|=\frac{1}{\|g_i\|^*}$, it follows that
$$\snorm[z_i]_{\widetilde{\Theta}}=\inf \{\|h\|_X \ : \ h\in X, 1\leq |g_i(h)|\}
=\frac{1}{\|g_i\|^*}.$$ 
Therefore, $$\|Vz_i\|=\|f_i\|\leq \frac{1}{\|g_i\|^*}= \snorm[z_i]_{\widetilde{\Theta}}, \ i\in\mn,$$ which implies that $V$ is bounded on $\seqgr[z]$.
Since $\widetilde{\Theta}$ is $CB$, extend $V$ on $\widetilde{\Theta}$ by linearity and continuity. For every $f\in X$ we have $\{g_i(f)\}_{i=1}^\infty\in\widetilde{\Theta}$ and $V(\{g_i(f)\}_{i=1}^\infty=V(\sumii g_i(f) z_i)=\sumii g_i(f) f_i =f$, which concludes the proof that $\seqgr[g]$ is a Banach frame for $X$ with respect to $\widetilde{\Theta}$. 
\ep

\v The sequences $\{g_i^1\}$ and $\{g_i^2\}$ from Example \ref{ex3}
give series expansion. These sequences are neither Hilbert frames, nor Hilbert Bessel sequences. Thus the existence of the series expansions do not follow neither from the Hilbert frame expansions, nor from \cite[Proposition 4.7]{PS} (which applies to Bessel sequences). Note that Proposition \ref{primer1new} applies to $\{g_i^1\}$ and $\{g_i^2\}$ - these two sequences are Banach frames with respect to the corresponding space $\widetilde{\Theta}$, given by (\ref{setm2}), and $\widetilde{\Theta}$ is $CB$.

\section{Construction of $\Theta_F$ and a class of $F$-frames}

In the next theorem we extend the construction of $\widetilde{\Theta}$ of Theorem \ref{constructone2} to the construction of $\Theta_F=\cap_{s\in\mn_0} \Theta_s$. The proof goes in the same way as in \cite[Theorem 4.8]{PS}, applying Theorem \ref{constructone2} to $X=X_s$ and
$\{g_i\vert_{X_s}\}_{i=1}^\infty$. 

\begin{thm} \label{constructthetas2}
  Let $\{X_s\}_{s\in\mn_0}$ be a sequence of reflexive Banach spaces which satisfies (\ref{fx1})-(\ref{fx3})
  and let $\seqgr[g]\newin (X^*)^\mn\setminus \{\nulel\}$. 
 For every $s\newin\mn_0$ and every scalar sequence $c=\seqgr[c]$, let $M^{c}_s$ be given by (\ref{mx2})
and define
 \begin{equation}\label{ts2new}
 {\Theta}_s := \left\{ c=\seqgr[c] \ : \ M^c_s\neq
\varnothing\right\}, \ \, \snorm[c]_{s} :=\inf\left\{ \|f\|_s : f\in M^c_s\right\}.
 \end{equation}
Assume that $(\mathcal{A}_1^s)$ holds for every $s\in\mn_0$. Then $(\mathcal{P}_1)-(\mathcal{P}_3)$ hold.
\end{thm}

Direct consequences of Theorem \ref{constructthetas2} are given in
 the next corollary.
\begin{cor}\label{c1} Let the assumptions of Theorem \ref{constructthetas2} hold. Then the following holds.

\noindent {\rm(a)} If $(\mathcal{A}_1^s)$ and $(\mathcal{A}_2^s)$ are satisfied for every $s\in\mn_0$, then $\{\Theta_s\}_{s\in\mn_0}$
is a sequence of solid $CB$-spaces with the properties (\ref{fx1})-(\ref{fx3}) such that $\{g_i|_{X_F}\}_{i=1}^\infty$ is an $F$-Bessel sequence
for $X_F$ with respect to $\Theta_F$.

\vspace{.1in}\noindent {\rm(b)} If $(\mathcal{A}_1^s)$, $(\mathcal{A}_2^s)$  and $(\mathcal{A}_3^s)$ are satisfied for every
$s\in\mn_0$, then $\{g_i|_{X_F}\}_{i=1}^\infty$ is a pre-$F$-frame for $X_F$ with respect to $\Theta_F$.
\end{cor}

The next proposition generalizes \cite[Proposition 4.10]{PS}. With this we construct a class of $F$-frames.

\begin{prop}\label{hilbcornew}
Let $(X_0, \<\cdot,\cdot\>_0)$ be a Hilbert space and let $\seqgr[e]$ denote an orthonormal basis for $X_0$. For given number sequences
$\{a_{i,s}\}_{i=1}^\infty$, $s\newin\mn$, with $1\leq a_{i,s}\leq a_{i,s+1}$, $i\newin\mn$, $s\newin\mn$, define
$$
X_s:=\left\{f\in X_0 \ : \ \{a_{i,s} \<f,e_i\>_0\}_{i=1}^\infty \in \ell^2\right\}, \ \ \<f,h\>_s:=\sumii a_{i,s}^2 \<f,e_i\>_0 \,\<e_i, h\>_0.
$$
Let $\Theta=\ell^2$ and $\seqgr[g]\newin (X_0^*)^\mn$ be defined by
\begin{eqnarray*}
g_i(f):=t_i\<f, e_1\>_0, &\ & i=1,2,\ldots, k_1, \\
g_i(f):=t_i\<f, e_j\>_0, &\ & i=k_{j-1}+1, k_{j-1}+2, \ldots, k_j, \ j\in\mn, j>1,
\end{eqnarray*}
where $k_j\in\mn$, $t_j\in\mr$, $t_j\neq 0$, $j\in\mn$.
 Then $\{X_s\}_{s\in\mn_0}$ is a sequence of Hilbert spaces, which satisfies
(\ref{fx1})-(\ref{fx3}); $\{\Theta_s\}_{s\in\mn_0}$, constructed by (\ref{ts2new}), is a sequence of $CB$-spaces, which satisfies
(\ref{fx1})-(\ref{fx3}) and $\{g_i|_{X_F}\}_{i=1}^\infty$ is a tight $F$-frame for $X_F$ with respect to $\Theta_F$.
\end{prop}
\bp In the same way as in \cite[Proposition 4.10]{PS}, it is not difficult to see that $\{X_s\}_{s\in\mn_0}$ is a sequence of Hilbert spaces,
which satisfies (\ref{fx1})-(\ref{fx3}) and such that for every $s\in\mn$, the sequence $\{e_i/a_{i,s}\}_{i=1}^\infty$ is an orthonormal basis for $X_s$. Denote $z_{i,s}:=e_i/a_{i,s}$, $i\in\mn$, $s\in\mn$

Let us now show that $(\mathcal{A}_1^s)$ is fulfilled. Take $c\newin\Theta_s$, $d\newin\Theta_s$, $f\newin M_s^c$, $h\newin M_s^d$.
Define $\widetilde{c}_i$, $\widetilde{d}_i$, $i\in\mn$, by (\ref{c1d1})-(\ref{dj}).
By (\ref{mx2}), 
$|c_i|\leq |g_i(f)| \ \mbox{and} \ |d_i|\leq |g_i(h)|,$ $i\in\mn.$ Therefore, $|\widetilde{c}_i|\leq |\<f,e_i\>_0|$, $i\newin \mn$, and 
$$ \sum_{i=1}^\infty a_{i,s}^2 \widetilde{c}_i^2
\leq\sum_{i=1}^\infty a_{i,s}^2|\<f,e_i\>_0|^2=\sum_{i=1}^\infty |\<f,z_i\>_s|^2= ||f||_s^2<\infty,
$$
which implies that
\begin{equation} \label{shod} \sum_{i=1}^\infty a_{i,s}^2\widetilde{c}_{i+1}^{\,\,2} < \infty;
\ \mbox{similarly,}\ \sum_{i=1}^\infty a_{i,s}^2\widetilde{d}_{i+1} ^{\,\, 2} < \infty.
\end{equation}

Let us find $r^{f,h}\in M_s^{c+d}$ such that $\|r^{f,h}\|_s\leq \|f\|_s+\|h\|_s$. Consider
$$r^{f,h}:= m_1 e_1 + m_2 e_2 + m_3 e_3 + m_4 e_4
+\ldots,$$ where $m_i$, $i\in\mn$, are given by (\ref{m1}) and (\ref{mj}).
Since $m_i\leq \widetilde{c}_i+\widetilde{d}_i$, (\ref{shod}) implies that $\{m_ia_{i,s} \}_{i=1}^\infty\in\ell^2$ and hence, the element
\begin{equation}\label{rfh}
r^{f,h}:=m_1a_{1,s}z_1+m_2a_{2,s}z_2 + m_3 a_{3,s}z_{3}+\ldots
\end{equation}
belongs to $X_s$.

 By definition, we have
\begin{equation}\label{econd}
 \<f,z_{i}\>_s= a_{i,s}\, \<f, e_i\>_0, \ f\in X_s,\ i\in\mn.
  \end{equation}
Now it is clear that for $i=1,2,\ldots, k_1$,
$$|c_i+d_i|\leq t_i\, m_1=
\frac{t_i}{a_{1,s}} \<r^{f,h},z_1\>_s=t_i\,\<r^{f,h},e_1\>_0 =g_i(r^{f,h}),$$
and for $i=k_{j-1}+1, k_{j-1}+2,\ldots, k_j,$  $j\in\mn$, $j>1,$
$$|c_i+d_i|\leq t_i\,m_j=\frac{t_i}{a_{j,s}} \<r^{f,h},z_j\>_s=t_i\,\<r^{f,h},e_j\>_0 =g_i(r^{f,h}).$$
This implies that $r^{f,h}\newin M^{c+d}$.
In a similar way as in Proposition \ref{primer1new} we obtain that
\begin{eqnarray*}
\|r^{f,h}\|_s &=&
\left\|\left\{m_1a_{1,s}, m_2a_{2,s}, m_3a_{3,s}, \ldots  \right\}\right \|_{\ell^2}\\
&\leq & \left\|\left\{a_{1,s}\widetilde{c}_1+a_{1,s}\widetilde{d}_1, a_{2,s}\widetilde{c}_2+a_{2,s}\widetilde{d}_2,
a_{3,s}\widetilde{c}_3+a_{3,s}\widetilde{d}_3,
\ldots \right\} \right\|_{\ell^2}\\
& \leq & \left\|\left\{a_{1,s}\widetilde{c}_1, a_{2,s}\widetilde{c}_2, a_{3,s}\widetilde{c}_3, \ldots
\right\}\right\|_{\ell^2} \\
 & &+ \left\|\left\{a_{1,s}\widetilde{d}_1,a_{2,s}\widetilde{d}_2,a_{3,s}\widetilde{d}_3, \ldots \right\}
\right\|_{\ell^2}\\
& \leq & \left\|\left\{a_{i,s}\<f,e_i\>_0\right\}_{i=1}^\infty\right\|_{\ell^2}  +\left\|\left\{a_{i,s}\<h,e_i\>_0 \right\}_{i=1}^\infty \right\|_{\ell^2}\\
& = & \left\|\left\{\<f,z_i\>_s\right\}_{i=1}^\infty\right\|_{\ell^2} +\left\|\left\{\<h,z_i\>_s \right\}_{i=1}^\infty
\right\|_{\ell^2}=\|f\|_s+\|h\|_s.
\end{eqnarray*}

Let us now prove that $(\mathcal{A}_3^s)$ is fulfilled. Consider $f\newin X_s$ and take arbitrary $\widetilde{f}\in
M_s^{\{g_i(f)\}_{i=1}^\infty}$, i.e. $|g_i(f)|\leq |g_i(\widetilde{f})|, i\in\mn$, and hence $|\<f,e_i\>_0|\leq |\< \widetilde{f},e_i\>_0|$,
$i\in\mn$. Using (\ref{econd}), we obtain
$$
\|f\|_s^2=\sumii a_{i,s}^2 \,|\<f, e_i\>_0|^2\leq \sumii a_{i,s}^2 \,|\<\widetilde{f}, e_i\>_0|^2= \|\widetilde{f}\|_s^2.
$$
Therefore $(\mathcal{A}_3^s)$ holds with $A_s=1$.

Take now $c=\{c_1,c_2,c_3,\ldots \}\in \Theta_s$. Fix $\varepsilon>0$ and by (\ref{shod}), find
$p\in\mn$, $p>1,$ such that $\sum_{i=k}^\infty a_{i,s}^2\widetilde{c}_i^{\,2}<\varepsilon$. Define
$$\seqgr[b]:=\{\underbrace{0, \ldots, 0}_{p-1}, a_{p,s}\widetilde{c}_p, a_{p+1,s}\widetilde{c}_{p+1},a_{p+2,s}\widetilde{c}_{p+2},\ldots \}.$$
Since $\seqgr[b]\in\ell^2$, there exists $h\newin X_s$ such that $b_i=\<h,z_i\>_s$, $i\in\mn$, and
$\|h\|_s^2=\sum_{i=1}^\infty |b_i|^2 = \sum_{i=p}^\infty a_{i,s}^2 \widetilde{c}_i^2< \varepsilon.$ Take $k:=k_{p-1}$ and recall,
$$c^{(k)}=\{\underbrace{0,\ldots,0}_{k_{p-1}},c_{k_{p-1}+1},c_{k_{p-1}+2},\ldots, c_{k_p},c_{k_p+1}, \ldots\}.$$
Thus, for $i=1,2,\ldots,k_{p-1}$,
$$|c^{(k)}_i|=0\leq |g_{i}(h)|$$ and for $i=k_{p+n-1}+1$, $\ldots, k_{p+n}$,  $n\in\mn_0$,
$$|c^{(k)}_i|\leq \widetilde{c}_{p+n}=\frac{b_{p+n}}{a_{p+1,s}}=\frac{\<h,z_{p+n}\>_s}{a_{p+n,s}} =\<h,e_{p+n}\>_0= g_{i}(h).
$$
This implies that $h\in M_s^{c^{(k)}}$. Since $\|h\|^2<\varepsilon $, it follows that $(\mathcal{A}_2^s)$ is fulfilled.
Now Corollary \ref{c1} implies that $\{g_i|_{X_F}\}_{i=1}^\infty$ is a tight pre-$F$-frame for $X_F$ with respect to $\Theta_F$. 

It remains to prove the existence of an $F$-frame operator. 
Denote the canonical vectors  by $z_i$, $i\in\mn$, and note that they 
form a basis for $\Theta_F$ in the sense that every $\seqgr[c]\in\Theta_F$ cab be written as $\seqgr[c]=\sumii c_i z_i$ with the convergence in $\Theta_s$-norm for every $s\in\mn_0$.
Let $\seqgr[f]$ and $\seqgr[h]$ be given by (\ref{seqf})-(\ref{hj}). Note that $f_i\in X_F$, $h_i\in X_F$, $i\in\mn$. 
Define $V$ on $\seqgr[z]$ by $Vz_i:=f_i$, $i\in\mn$. Our aim is to prove that for every $s\in\mn_0$, $V$ is $s$-bounded on $\seqgr[z]$.
Fix $s\in\mn_0$. Let us first prove that $\snorm[z_i]_s=1/\|g_i\vert_{X_s}\|^*_s$, $i\in\mn$.
  Fix $i\in\mn$ and note that $h\in M^{z_i}_s$ if and only if $1\leq |g_i(h)|$. Thus, the elements $h\in M^{z_i}_s$ should satisfy the inequality $\|h\|_s\geq \frac{1}{\|g_i\vert_{X_s}\|^*_s}$. 
Since $h_i\in M^{z_i}_s$ and $\|h_i\|_s=\frac{1}{\|g_i\vert_{X_s}\|^*_s}$, it follows that
$$\snorm[z_i]_s=\inf \{\|h\|_s \ : \ h\in X_s, 1\leq |g_i(h)|\}
=\frac{1}{\|g_i\vert_{X_s}\|^*_s}.$$ 
Therefore, $$\|Vz_i\|_s=\|f_i\|_s\leq \frac{1}{\|g_i\vert_{X_s}\|^*_s}= \snorm[z_i]_s, \ i\in\mn,$$ which implies that $V$ is $s$-bounded on $\seqgr[z]$.
Let $V_s:\Theta_s\to X_s$ denote the extension of $V$ on $\Theta_s$ by linearity and continuity.
Note that $V_0\vert_{\Theta_F}=V_s\vert_{\Theta_F}$, $s\in\mn$. For every $c=\seqgr[c]\in\Theta_F$, the sequence $\{\sum_{i=1}^n c_i f_i\}_{n=1}^\infty$ converges in $X_s$-norm (to $V_s c$)  for every $s\in\mn_0$. Thus,  $\{\sum_{i=1}^n c_i f_i\}_{n=1}^\infty$ converges in $X_F$. Therefore, $V_0\vert_{\Theta_F}$ is an $F$-bounded operator from $\Theta_F$ into $X_F$. 
Let $s\in\mn_0$.
For every $f\in X_F$, $\sum_{i=1}^n g_i(f) f_i\to f$ in $\|\cdot\|_s$-norm 
and $\sum_{i=1}^n g_i(f) f_i=V_s(\sum_{i=1}^n g_i(f) z_i)\to V_s(\{g_i(f)\}_{i=1}^\infty)$ in $\|\cdot\|_s$-norm as $n\to\infty$. Since this holds for every $s\in\mn$,
we have that $V_0\vert_{\Theta_F}(\{g_i(f)\}_{i=1}^\infty)=f$, $f\in X_F$. This concludes the proof that $V_0\vert_{\Theta_F}$ is an $F$-frame operator for $\seqgr[g]$.
\ep

\end{document}